\documentclass[11pt]{amsart}
\usepackage{amssymb,latexsym,amsfonts}
\usepackage{amsmath}
\usepackage{amscd}
\usepackage{pb-diagram,pb-xy}%lamsarrow,pb-lams}

\usepackage{graphicx}
\usepackage{hyperref}
\usepackage{color}
\usepackage[matrix,arrow]{xy}

\usepackage{fullpage}
\newcommand{\vol}{{\rm Vol}}

\newcommand{\Ein}{{\mathrm {Ein}}}
\newcommand{\cEin}{{\mathbf {Ein}}}
\newcommand{\ein}{{\mathrm {ein}}}
\newcommand{\cein}{{{\mathbf {ein}}}}
\newcommand{\Eink}{{\mathrm {Ein}_k}}

\newcommand{\Einl}{{\mathrm {Ein}_l}}

\newcommand{\Scal}{{\rm Scal}}
\newcommand{\tr}{\operatorname{tr}}

\newcommand{\Ric}{{\rm Ric}}

\DeclareMathAlphabet{\mathpzc}{OT1}{pzc}{m}{it}

\def\rhom{\rho_{{\rm max}}}
\def\rhomin{\rho_{{\rm min}}}

\newtheorem{theorem}{Theorem}[section]

\newtheorem{proposition}[theorem]{Proposition}
\newtheorem*{thm-a}{Theorem\! A}
\newtheorem*{thm-aa}{Theorem\! ${\rm A}^{\prime}$}
\newtheorem*{cor-bb}{Corollary\! ${\rm B}^{\prime}$}
\newtheorem*{thm-cc}{Theorem\! ${\rm C}^{\prime}$}
\newtheorem*{cor-a}{Corollary\! A}
\newtheorem*{thm-b}{Theorem\! B}
\newtheorem*{cor-b}{Corollary\! B}
\newtheorem*{thm-c}{Theorem\! C}
\newtheorem*{cor-c}{Corollary\! C}
\newtheorem*{thm-d}{Theorem\! D}
\newtheorem*{thm-dd}{Theorem\! ${\rm D}^{\prime}$}
\newtheorem*{cor-dd}{Corollary\! ${\rm D}^{\prime}$}
\newtheorem*{cor-d}{Corollary\! D}
\newtheorem*{thm-e}{Theorem\! E}
\newtheorem*{cor-e}{Corollary\! E}
\newtheorem*{thm-f}{Theorem\! F}
\newtheorem*{cor-f}{Corollary\! F}

\newtheorem*{conj-d}{Conjecture D}
\newtheorem*{conj-c}{Conjecture C}
\newtheorem{thm}{Theorem}[section]
\theoremstyle{definition}
\newtheorem{definition}[thm]{Definition}

\newtheorem*{remark}{Remark}
\newtheorem*{remarks}{Remarks}

\newtheorem{example}{Example}[section]

\begin{document}
\author{Mohammed Larbi Labbi} 
\address{Department of Mathematics\\
 College of Science\\
  University of Bahrain\\
  32038, Bahrain.}
\email{mlabbi@uob.edu.bh}
%\subjclass[2020]{53C21, 53C18}
\renewcommand{\subjclassname}{
  \textup{2020} Mathematics Subject Classification}

\title{On the conformal Ein invariants}  \maketitle
 \begin{abstract}
For a compact  Riemannian $n$-manifold $(M,g)$ of positive scalar curvature, the capital $\Ein$ invariant of $g$ is defined to be the infinimum over $M$ of the quotient of the scalar curvature by the maximal eigenvalue of the Ricci curvature. This is a re-scale invariant and belongs to the interval $(0,n]$. For a positive conformal class $[g]$, we define the conformal invariant $\Ein([g]):=\sup\{\Ein(g):  g\in [g]\}$. In this paper, we prove vanihing theorems for  Betti numbers and for the higher homotopy groups of $M$ under  optimal lower bounds on $\Ein([g])$ assuming that $g$ is locally conformally flat. We establish an inequality relating our invariant to Schoen-Yau conformal invariant $d(M,[g])$  from which we deduce a classification result for locally conformally flat manifolds with higher $\Ein([g])$. We show that  the class of locally conformally flat manifolds with $\Ein([g])>k$ is stable under the operation of connected sums for  $0<k<n-1.$\\
For a general positive conformal class, we prove in dimension $4$ an inequality relating $\Ein([g])$ to the first and second Yamabe invariants.
 Similar results are proved in this paper for an analogous conformal invariant, namely the small $\ein$ invariant.
\end{abstract}

\tableofcontents
\section{Introduction }

Throughout this paper, $(M,g)$ denotes a closed connected Riemannian manifold of dimension $n$. We denote   by $\Ric$ and $\Scal$ the Ricci  and  scalar curvatures of $(M,g)$. Recall that the $k$-th modified Einstein tensor on $(M,g)$ is defined by 
\begin{equation}
\Eink (g):=\Scal \, g -k\Ric,
\end{equation}
Where $k$ is a constant. We bring to the attention of the reader that this modification was used first by Rastall in 1972 in his modified theory of Gravity \cite{Rastall} and studied by Bourguignon in \cite{Bourguignon}.\\
We are interested in  the positivity properties of these modified tensors. We first remark that for a compact Riemannian $n$-manifold $(M,g)$ with positive scalar curvature and for $0<k<n$, the tensor $\Eink (g)$ is positive definite if and only if at each point of $M$ one has
\begin{equation}\label{rhom}
k<\frac{\Scal}{\rhom},
\end{equation}
 Where $\rhom$ denotes the maximal eigenvalue of the Ricci curvature. 
\subsection{The capital $\Ein$ invariant}
A straightforward consequence of the above charecterisation of the positivity of the  tensors $\Eink (g)$ is the following descent positivity property
\begin{equation}
{\rm For}\,\, 0<k<l<n,\,\,\, \Einl>0 \Rightarrow \Eink>0.
\end{equation}

We therefore define the metric invariant 
\begin{equation}
\Ein(g):=\sup\{k \in (0, n): \Eink(g)>0\}.
\end{equation}
We set $\Ein(g)=0$  if the scalar curvature of $g$ is not positive. An immediate consequence of \ref{rhom}, is that for a metric $g$ of positive scalar curvature one has
\begin{equation}
\Ein(g)=\inf_{M}\frac{\Scal(g)}{\rho_{{\rm max}}(g)}.
\end{equation}
The above formula shows in particular that the metric invariant $\Ein(g)$ is re-scale invariant, That is for any positive real number $t$, one has
\begin{equation}
\Ein(tg)=\Ein(g).
\end{equation}
Recall for later reference that \cite{modified-Eink}
\begin{equation}
\begin{split}
\Ein(g)>n-1\implies \Ric >0.\\
\Ein(g)=n-1\implies \Ric \geq 0.
\end{split}
\end{equation}
Conversely, $\Ric>0$ implies only that $\Ein >1$ and this is optimal as shown by the example of Berger metrics on $S^{2n+1}$.
The ${\rm Ein}$ invariant defines a {\emph{pre-order}} on the set of Riemannian metrics on $M$:
\begin{equation}
g_1\preceq g_2 \,\,\, {\rm if}\,\,\, \Ein(g_1)\leq \Ein(g_2).
\end{equation}
This lead us naturally to the study of maximal metrics with respect to the above pre-order. \\
We define the smooth invariant $\Ein(M)$ to be
\begin{equation}
\cEin(M)= \sup\{\Ein(g)\colon  g\in {\mathcal M}\},
\end{equation}
where ${\mathcal M}$ denotes the space of all  Riemannian metrics on $M$.

\subsection{The small $\ein$ invariant}
The tensors $\Eink (g)$ for negative $k$ enjoy the  following descent positivity property
\begin{equation}
{\rm For}\,\, l<k<0,\,\,\, \Einl>0 \Rightarrow \Eink>0\Rightarrow \Scal >0.
\end{equation}
We then naturally define \cite{modified-Eink} the metric  invariant $\ein(g)$ to be 
\begin{equation}
\ein(g):=\inf\{k<0\colon \Eink(g)>0\}
\end{equation}
We set $\ein(g)=0$ if the scalar curvature of $g$ is not positive and $\ein(g)=-\infty$ in case the corresponding set of $k$'s is unbounded below. \\
We remark that if the metric $g$ has nonnegative Ricci curvature and positive scalar curvature   then clearly $\ein(g)=-\infty$. Otherwise, $\ein(g)$ can be alternatively defined by
\begin{equation}
\ein(g)=\sup_{x\in M}\Bigl\{\frac{\Scal(x)}{\rhomin(x)}:\rhomin(x)<0\Bigr\},
\end{equation}
where $\rhomin(x)$ denotes the minimum eigenvalue of Ricci curvature of $g$ at $x\in M$.
\begin{remark}
We bring to the attention of the reader that Guan and Wang \cite{Guan-Wang}, proved that in the case $\rhomin$ is everywhere negative on $M$, the existence of a unique conformal metric $\bar{g}$ whose minimum of Ricci curvature is constant over $M$ and  equal to $-1$. In particular,  if the scalar curvature of $\bar{g}$ is positive then $\ein(\bar{g})$ is nothing but the negative of the infinimum over $M$ of the scalar curvature of $\bar{g}$.
\end{remark}
A second smooth invariant $\cein(M)$  of $M$ is defined as follows
\begin{equation}
\cein(M)= \inf \{\ein(g)\colon  g\in {\mathcal M}\}.
\end{equation}
The above  invariants were studied in  \cite{modified-Eink}. In this paper, we restrict our study of the above constants once restricted to a given conformal class.   We start by a definition

\begin{definition}
Let $[g]$ denotes the  conformal class of the metric $g$, we define the following two constants of $[g]$
\begin{equation}
\Ein([g])= \sup\{\Ein(g):  g\in [g]\} \quad {\rm and}\quad \ein([g])= \inf \{\ein(g):  g\in [g]\}.
\end{equation}
\end{definition}

We remark first that if the conformal class of $g$ contains an Einstein metric with positive scalar curvature then $\Ein([g])=n$. Similarly, if the conformal class of $g$ contains a metric with non-negative Ricci curvature and with positive scalar curvature then $\cein([g])=-\infty$.
\begin{definition}\label{def1}
\begin{itemize}
\item A metric $g$ is said to be $\Ein(M)$-maximal if $\Ein(g)=\Ein(M).$
\item A metric $g$ is said to be $\Ein([g])$-maximal if $\Ein(g)=\Ein(M).$
\end{itemize}
\end{definition}
\section{Statement of the main results}
The first result is a vanishing theorem
\begin{thm-a}\label{vanish-theorem}
Let $(M,g)$ be a compact locally conformally flat $n$-manifold  and $p$ an integer such that $1\leq p\leq \frac{n}{2}$. 
\begin{enumerate}
\item  If $\Ein([g])>\frac{(n-1)(n-2p)}{n-p-1}$  then the Betti numbers  $b_k$ of $M$ vanish for $p\leq k\leq n-p$.
\item If $p>1$ and  $\ein([g])<\frac{(n-1)(2p-n)}{p-1}$ then the Betti numbers  $b_k$ of $M$ vanish for $p\leq k\leq n-p$.
\end{enumerate} 
\end{thm-a}
As a consequence we are able to compute the $\Ein([g])$ and $\ein([g])$ constants for the product of two space forms of opposite signs as follows
\begin{cor-a}\label{SnHd}
Let $n,d$ be positive integers such that $d<\frac{n-2}{2}$. Let  $(S^{n-d-1},g_0)$ be the standard sphere of curvature
$+1$ and let $(M^{d+1},g_1)$ be a compact space form of curvature $-1$.
We denote by $g$ the Riemannian product of $g_0$ and $g_1$ on $S^{n-d-1}
\times M^{d+1}$. Then 
\[\Ein([g])=\Ein(g)=(n-1)\frac{2d-n+2}{d-n+2}\quad {\rm and}\quad \ein([g])=\ein(g)=-(n-1)\frac{n-2d-2}{d}.\]
In other words the metric $g$ is $\Ein([g])$-Maximal and $\ein([g])$-Minimal.
\end{cor-a}
The next theorem provides a converse result to the previous corollary and the limit case in the above theorem
\begin{thm-aa}
Let $(M,g)$ be a compact  locally conformally flat $n$-manifold with $n>2$,  and $p$ an integer such that $1\leq p< \frac{n}{2}$. 
\begin{enumerate}
\item  If $\Ein(g)=\frac{(n-1)(n-2p)}{n-p-1}$ and the Betti number  $b_p\not= 0$ then $M$ is covered by the standard  product $S^{n-p}\times H^p$ of the $(n-p)$-sphere and the hyperbolic space of dimension $p$.
\item If $p>1$, $\ein(g)=\frac{(n-1)(2p-n)}{p-1}$ and the Betti number  $b_p\not= 0$ then $M$ is covered by the  standard product $S^{n-p}\times H^p$.
\item If $p=1$, $\ein(g)=-\infty $ and the Betti number  $b_1\not= 0$ then $M$ is covered by the standard product $S^{n-1}\times \mathbb{R}$.
\end{enumerate} 
\end{thm-aa}
Let now $(M,g)$ be a smooth compact connected  locally conformally flat Riemannian $n$-manifold with positive scalar curvature of dimension $n\geq 4$, we denote by $\tilde{M}$ its universal cover. Schoen and Yau \cite{Schoen-Yau}, see also \cite{Izeki, Izeki2, Nayatani}, proved then that the developing map $\Phi: \tilde{M} \to S^n$ is a conformal embedding, $\pi_1(M)$ is isomorphic to a discrete subgroup $\Gamma$ of ${\rm Conf}(M)$, $\Phi(\tilde{M})$ is a domain $\Omega$ in $S^n$ and coincides with the complement $\Omega(\Gamma)$ of the limit set $\Lambda(\Gamma)$ of the action of $\pi_1(M)\approx \Gamma$ on the sphere $S^n$. In other words, $(M,g)$  is conformally equivalent to the Kleinian manifold $\Omega/\Gamma$.\\
Let $\delta:=\delta(\Gamma)$ denotes the critical exponent of the Kleinian group $\Gamma$, see \cite{Schoen-Yau,Nayatani, Izeki, Izeki2}. It turns out that $\delta$ depends only on the conformal class of $g$, precisely $\delta$ coincides with the Schoen-Yau conformal invariant $d(M, [g])$ of the conformal class $[g]$ of $g$, see Theorem 4 in \cite{Izeki2}. If in  addition the Kleinian group $\Gamma$ is not elementary then $\delta$ coincides with the Hausdorff dimension of the limit set $\Lambda(\Gamma)$.\\
Shoen-Yau proved in \cite{Schoen-Yau} that for the above manifold one has $d(M, [g])\leq \frac{n-2}{2}$. In the next theorem we generalize this result as follows 
\begin{thm-b}
Let $(M,g)$ be a compact connected locally conformally flat Riemannian manifold of dimension $n\geq 4$ such that $\Ein([g])>0$. Then 
$$d\leq (n-2)\frac{n-1- \Ein([g])}{2(n-1)- \Ein([g])}.$$
Where $d$ denotes the Schoen-Yau conformal invariant $d(M, [g])$ of the conformal class   $[g]$. Equivalently, we have
$$ \Ein([g])\geq (n-1)\frac{2d-n+2}{d -n+2}.$$
\end{thm-b}
Recall that the condition $\Ein([g])>0$ is equivalent to the existence of a conformal metric in $[g]$ of positive scalar curvature.\\
The above inequality is optimal as shown by the standard metric on the product $S^{n-d-1}\times \mathbb{H}^{d+1}$. For more details see Example \ref{Example1}.

\begin{cor-b}
Let $(M,g)$ be a compact connected locally conformally flat Riemannian manifold  of dimension $n\geq 4$. Let $p$ be an integer such that $2\leq p\leq n-1$. Then one has 
\begin{equation*}
\Ein[g]>\max\bigl\{0, \frac{(n-1)(2p-n)}{p-1}\bigr\}\implies \pi_2(M)=\pi_3(M)=...=\pi_p(M)=0.
\end{equation*}
\end{cor-b}

The next Corollary provides a classification of  locally conformally flat Riemannian manifold  of dimension $n\geq 4$ and of higher $\Ein$
\begin{cor-bb}
Let $(M,g)$ be a compact connected locally conformally flat Riemannian manifold  of dimension $n\geq 4$. If 
$\Ein[g]>\frac{(n-1)(n-4)}{n-3}$ then there is a finite covering of $M$ which is either diffeomorphic to $S^n$ or to connected sums of copies of $S^1\times S^{n-1}$.\\
The above inequality is optimal as shown by the example of the standard product $S^{n-2}\times \mathbb{H}^{2}$ whose conformal $\Ein$ 
equals $\frac{(n-1)(n-4)}{n-3}$.
\end{cor-bb}

\begin{remark}
It is remarkable that the  lower bound $\frac{(n-1)(n-4)}{n-3}$ for $\Ein([g])$ in Corollary B' is also the optimal lower bound required for $\Ein(g)$ in order for it to imply positive isotropic curvature (PIC). Precisely, it is not difficult to check that
for a locally conformally flat Riemannian manifold  of dimension $n\geq 4$ one has
$$\Ein(g)>\frac{(n-1)(n-4)}{n-3}\implies {\rm PIC}.$$ 
The same conclusion does not hold for $\Ein[g]\leq \frac{(n-1)(n-4)}{n-3}$ as shown by the example of the standard product $S^{n-2}\times \mathbb{H}^{2}$.
\end{remark}

The next Theorem C shows that within the class of conformally flat manifolds,  the positivity of modified Einstein tensors is preserved under connected sums
\begin{thm-c}\label{csum}
\begin{enumerate}
\item
For $k\in (-\infty, 2)$, the connected sum of two conformally flat manifolds each one of positive modified Einstein tensor $\Eink$  and with dimension $\geq 3$ admits a conformally flat metric  of positive $\Eink$ tensor.
\item For $k\geq 2$, the connected sum of two conformally flat manifolds each one of positive $\Eink$ tensor and with dimension $> k+1$ admits a conformally flat metric  of positive $\Eink$ tensor.
\end{enumerate}
\end{thm-c}
As an application of the above theorem one has,

\begin{cor-c} Let $n\geq 3$ and $p\geq 1$. For each real numbers $k_1\in [0, n-1), k_2\in (-\infty, 0]$, the connected sum of
 $p$ copies of $S^1\times S^{n-1}$ admits conformally flat metrics $g_1$ and $g_2$ such that $\Ein(g_1)>k_1$ and $\ein(g_2)<k_2$.
\end{cor-c}
\begin{remark}
For $p>1$,  the connected sum of $p$ copies of $S^1\times S^{n-1}$ can't admit a conformally flat metric $g$ with $\Ein(g)=n-1$ as this will imply that $g$ has non-negative Ricci curvature which is impossible, see \cite{Noronha}. It would be interesting to decide whether there exists on the previous connected sum a conformally flat metric $g$  with $\Ein([g])=n-1$.
\end{remark}
Next, we discuss the case of non conformally flat classes in four dimensions.\\
Let $(M,g)$ be a compact Riemannian $4$-manifold and let $A$ denotes its Schouten tensor, we denote as usual by $\sigma_1(A)$ and $\sigma_2(A)$ respectively the trace of $A$ and second elementary symmetric function in the eigenvalues of $A$. We recall two important conformal invariants of $[g]$. The first one is the celebrated Yamabe invariant $Y[g]$ and is defined by
\[Y[g]=\inf_{\bar{g}\in[g]}\frac{1}{\left(\vol(\bar{g})\right)^{1/2}}\int_M\Scal({\bar{g}})\mu_{\bar{g}}.\]
The second one is $\int_M\!\sigma_2(A)\mu_g$, that is  the integral over $M$ of $\sigma_2(A)$.\\
The positivity of  $\sigma_1(A)$ and $\sigma_2(A)$ together imply simultaneously the positivity of the Einstein tensor and the positivity of the Ricci tensor \cite{CGY}. Consequently they imply $\Ein([g])>2$ and $\ein([g])=-\infty$. The previous simple algebraic property was generalized in \cite{CGY} and \cite{Gursky-Viacklovsky} to  conformally invariant properties of the conformal class of $g$, see the main Theorem of \cite{Gursky-Viacklovsky} and Theorem A in \cite{CGY}.  As a consequence of their results we prove the following

\begin{thm-d}
Let $(M, g)$ be a compact $4$-dimensional Riemannian
manifold with positive  Yamabe invariant $Y[g]$.
\begin{enumerate}
\item If $\int\!\!\sigma_2(A) \mu_g>0 $ then $\Ein([g])>2$ and $\ein([g])=-\infty$.
\item If $\int\!\!\sigma_2(A)\mu_g =0$ then $\Ein([g])\geq 2$ and $\ein([g])=-\infty$.
\item If $\int\!\!\sigma_2(A)\mu_g<0$ then 
$$\Ein([g])\geq \frac{4Y[g]}{Y[g]+\sqrt{\left(Y[g]\right)^2-96\int\!\!\sigma_2(A)\mu_g}}\quad {\rm and}\quad  \ein([g])\leq \frac{4Y[g]}{Y[g]-\sqrt{\left(Y[g]\right)^2-96\int\!\!\sigma_2(A)\mu_g}}.$$
\end{enumerate}
\end{thm-d}

For a Riemannian manifold $(M,g)$ of dimension $4$, recall that the Paneitz operator $P_g$ is a fourth order generalization of the usual  Laplacian $\Delta$  defined  by
\begin{equation}
P_g(\phi)=\Delta^2\phi+\frac{2}{3} \delta\bigl({\rm Ein}_3(g) \bigr)d\phi.
\end{equation}
Where ${\rm Ein}_3(g)=\Scal \, g-3\Ric$.
An important future of $P_g$ is that it is conformally invariant. Precisely, if $\bar{g}=e^{-2u}g$ then
\begin{equation}
P_{\bar{g}}=e^{4u}P_g.
\end{equation}
In particular, the positivity of the Paneitz operator is a conformally invariant property and its kernel is conformally  invariant as well. The following theorem is a slightly weaker form of a Theorem  due to Eastwood and Singer,  see Theorem 5.5 in \cite{ES}, and also it is due to Gursky and Viacklovsky, see Proposition 6.1 in \cite{Gursky-Viacklovsky}. It provides a sufficient condition for the non-negativity of Paneitz operator.\\
\begin{thm-dd}
If a compact Riemannian $4$-manifold $(M,g)$ has  the property   $\Ein([g])>1$  then its Paneitz operator $P_g$ is non negative. Furthermore, the kernel of $P_g$ consists only of constant functions.
\end{thm-dd}
Associated with the Paneitz operator is the $Q$-curvature defined by
\begin{equation}
Q_g=-\frac{1}{12}\Delta {\Scal}_g+2\sigma_2(A_g).
\end{equation}
One consequence of the previous theorem is
\begin{cor-dd}
If a compact Riemannian $4$-manifold $(M,g)$ has  the property   $\Ein([g])>1$
 then it admits a conformal metric with constant $Q$-curvature.
\end{cor-dd}

\section{$\Ein([g])$ and $\ein([g])$ constants of a locally conformally flat class}
\subsection{Proof of Theorem A and Corollary A}

We  first prove Theorem A as follows
\begin{proof}
The hypothesis of the theorem implies the existence of a metric $g$ in the conformal class $[g]$ such that $\Eink(g)>0$ where $k=\frac{(n-1)(n-2p)}{n-p-1}$ in part (1), and $k=\frac{(n-1)(2p-n)}{p-1}$ in part (2).  To prove the theorem we use the descent positivity properties of the modified Einstein tensors and the classical Weitzenb\"ock formula for differential forms
$$\Delta = {\nabla}^*\nabla+{\mathcal {W}},$$
where  $\Delta$ is  the Laplacian of differential forms, $\nabla$ is the Levi-Civita connexion and ${\mathcal W}$ is the Weitzenb\"ock curvature term. We shall prove that under the theorem hypotheses the Weitzenb\"ock curvature term is positive and the theorem follows from the previous formula. The curvature term $ {\mathcal W}$ once operating on $p$-forms takes the following form \cite{Labbi-nachr}
\begin{equation}
{\mathcal W}_p=\frac{g^{p-2}}{(p-2)!}\left ( \frac{g\Ric}{p-1}-2R\right ).
\end{equation}
Where $R$ is the Riemann tensor and all the products are exterior products of double forms. Since the manifold is locally conformally flat then the Weyl tensor vanishes and we have $R=gA$ where
$$A=\frac{1}{n-2}\left ( \Ric -\frac{\Scal}{2(n-1)}g\right ),$$
is the Schouten tensor. Consequently, we can see that
\begin{equation}
\begin{split}
{\mathcal {W}}_p=& \frac{g^{p-2}}{(p-2)!}\left ( \frac{g\Ric}{p-1}-2gA\right )=
 \frac{g^{p-1}}{(p-1)!}\left( \Ric-(2p-2)A\right )\\
&= \frac{g^{p-1}}{(p-1)!}\left ( (1-\frac{2p-2}{n-2})\Ric+\frac{(p-1)\Scal}{(n-1)(n-2)}g\right)\\
&=\frac{p-1}{(n-1)(n-2)} \frac{g^{p-1}}{(p-1)!}\left ( \Scal\,  g -\frac{(n-1)(2p-n)}{p-1}\Ric \right )\\
&=\frac{p-1}{(n-1)(n-2)} \frac{g^{p-1}}{(p-1)!}\left ({\rm Ein}_{k} \right )\\
\end{split}
\end{equation}
Note that the last term  $\frac{g^{p-1}}{(p-1)!}\left (\Eink \right )$  is positive if and only if the sum of the lowest $p$ eigenvalues of $\Eink$ is positive. In particular, $\frac{g^{p-1}}{(p-1)!}\left (\Eink \right )$ is positive if $\Eink$ is positive. This proves the part (2) of the theorem. To prove the first part, it suffices to notice that the curvature term $ {\mathcal W}$ once operating on $(n-p)$-forms takes the form
\[ {\mathcal {W}}_{n-p}=\frac{n-p-1}{(n-1)(n-2)} \frac{g^{n-p-1}}{(n-p-1)!}\left ( \Scal\,  g -\frac{(n-1)(n-2p)}{n-p-1}\Ric \right ).\]
Part (2) follows then using Poincar\'e duality for $p\not=1$. If $p=1$, the condition $\Ein([g])>n-1$ implies the  the existence of a metric $g$ in the conformal class $[g]$  of positive Ricci curvature \cite{modified-Eink}. This completes the proof of Theorem A.
\end{proof}

Next, we prove Corollary A.

\begin{proof}
 A straightforward computation shows that for the standard product metric $g$ one has  $\Ein(g)=(n-1)\frac{2d-n+2}{d-n+2}$ and $\ein(g)=-(n-1)\frac{n-2d-2}{d}$. Therefore $\Ein([g])\geq(n-1)\frac{2d-n+2}{d-n+2}$ and $\ein([g])\leq-(n-1)\frac{n-2d-2}{d}$.\\
Next we use the notations of the theorem, let $p=n-d-1$ then $2p-n=n-2d-2>0$ and $k_1=\frac{(n-1)(2p-n)}{p-1}=(n-1)\frac{2d-n+2}{d-n+2}$ and 
$k_2=\frac{(n-1)(n-2p)}{n-p-1}=-(n-1)\frac{n-2d-2}{d}$.\\
The metric $g$ being conformally flat then any metric in the
conformal class of $g$ is conformally flat as well. Now since $b_{(d+1)}\left ( S^{n-d-1}
\times M^{d+1}\right )=b_{(d+1)}\left ( M^{d+1}\right )\not= 0$ then by Theorem A no conformally flat metric on $S^{n-d-1}
\times M^{d+1}$ can have $\Ein([g]) >k_1$ or $\ein([g])<k_2.$ This completes the proof.
\end{proof}

\subsection{Proof of Theorem ${\rm A}^{\prime}$}
\begin{proof}
We prove the capital $Ein$ part of the theorem. The small $\ein$ part is completely similar and left to the interested reader.\\
We proceed as in \cite{Lafontaine}, see also \cite{Noronha} and \cite{Labbi-Betti}. We have two possibilities. If $(M,g)$ is locally reducible, then it is either flat or covered by the standard product $S^k\times H^q$ \cite{Lafontaine}. Since $\Ein(S^k\times H^k)=\frac{(n-1)(n-2q)}{n-q-1}$ where $n=k+q$. With our hypothesis we must have $k=n-p$ and $q=p$.\\
If $(M,g)$ is locally irreducible, we distinguish two cases:\\
If the holonomy group is not $SO(n)$ or $U(n/2)$, then it follows from Berger's classification of holonomy groups that $(M,g)$ is Einstein, see for instance \cite{Besse}. Since $\Ein >0$ then it must be equal to $n$ which is not our case.\\
The second case is when the holonomy group is $SO(n)$ or $U(n/2)$. The hypothesis $\Ein(g)=\frac{(n-1)(n-2p)}{n-p-1}$ implies that the Weitzenb\"ock curvature operator of order $p$ is nonegative, so the Weitzenb\"ock formula shows that $M$ has a parallel $p$-form, hence  invariant under the action of the holonomy group. Therefore, at each point of the manifold, we have a $p$-form invariant under the action of $SO(n)$ or $U(n/2)$. The first possibility can't occur as $\wedge^p \mathbb{R}^n$ has no invariant subspaces of dimension $1$ under the action of $SO(n)$. In the second situation, $(M,g)$ must be Kahlerian and conformally flat and consequently flat if $n>4$, see 2.68 in \cite{Besse}. Remains to prove the theorem in the case $n=4$. In this case $p=1$ and $\Ein(g)=3=n-1$, consequently the Ricci curvature of $(M,g)$ is nonnegative and the theorem follows from the classification of conformally flat manifolds with nonnegative Ricci curvature \cite{Noronha}.
\end{proof}

\subsection{Nayatani metric and the $\Ein([g])$ invariant: Proof of Theorem B and Corollaries B, ${\rm B}^{\prime}$}
\subsubsection{Proof of theorem B}
\begin{proof}
It follows from the  discussion preceding the statement of Theorem B that $M=\Omega/\Gamma$ is a Kleinian manifold. We suppose that  $\delta(\Gamma)>0$. Nayatani \cite{Nayatani} constructed a canonical conformally flat metric $\bar{g}\in[g]$ on $M=\Omega/\Gamma$ whose Ricci curvature is  given by
\begin{equation}
\Ric = -(n - 2)(\delta + 1){\mathcal A} + (n -2-\delta)(\tr_{\bar{g}}{\mathcal A}) \bar{g},
\end{equation}
where ${\mathcal A}$ is a non-negative tensor whose trace $\tr_{\bar{g}}{\mathcal A}$ is strictly positive as we supposed the scalar curvature is positive, see \cite{Nayatani}.

Consequently, one has
\begin{equation}
\Eink(\bar{g})=\left((n-1)(n-2-2\delta)-k(n-2-\delta)\right)(\tr_{\bar{g}}{\mathcal A}) \bar{g}+k(n - 2)(\delta + 1){\mathcal A}.
\end{equation}
This is clearly positive if $k<(n-1)\frac{2\delta-n+2}{\delta -n+2}.$\\
In case $\delta=0$, then $(M,g)$ admits a conformal metric $\bar{g}$ with which $M$ is locally isometric to $S^{n-1}\times {\Bbb R}$ \cite{Nayatani}. Therefore, in this case one has $\Ein([g])\geq \Ein(\bar{g})=n-1$. Finally, If $\delta=-1$, that is $\pi_1(M)$ is finite, one has  $\Ein([g])=n$.

\end{proof}
The following example shows the optimality of the inequality of Theorem B.
\begin{example}\label{Example1}
Let $S^d$ be a round $d$-sphere in the $n$-sphere $S^n$ with $1\leq d\leq n-3$. Let $\Omega=S^n\backslash S^d$ and $G={\rm Conf}(\Omega):=\{f\in
{\rm Conf}(S^n): f(S^d)=S^d\}.$ It turns out that there exists a $G$-invariant conformally flat metric, say $g_0$, on $\Omega$ that makes it isometric to the standard product $S^{n-d-1}\times H^{d+1}$, see \cite{Nayatani}.\\
Let now $\Gamma\subset G$ be a Kleinian group with limit set   $\Lambda(\Gamma)=S^d$. The manifold $M=\Omega/\Gamma$ is conformally flat and here $\delta=\delta(\Gamma)$ coincides with the Hausdorff dimension $d$ of the limit set $S^d$. The so obtained metric on $M$ coincides with Nayatani metric up to a positive constant factor. Corollary  \ref{SnHd} shows that 
$$ \Ein([g_0])=\Ein(g_0)= (n-1)\frac{2d-n+2}{d -n+2}.$$

\end{example}

\subsubsection{Proof of Corollary B}
\begin{proof}
Let $(M,g)$ be compact and locally conformally flat $n$-manifold, and let $p$ be an integer such that $2\leq p\leq n-1$. Denote by $d$ the Schoen-Yau conformal invariant of $(M,[g])$.\\
We distinguish two cases. If $n\geq 2p$, we have $\Ein[g]>0$ and therefore $[g]$ contains a locally conformally flat metric of positive scalar curvature. The corollary follows then from Schoen-Yau theorem, namely Theorem 4.6(ii) in  \cite{Schoen-Yau}. In the case where $n<2p$, one has $\Ein[g]>\frac{(n-1)(2p-n)}{p-1}$ and this implies that
$$(n-2)\frac{n-1- \Ein([g])}{2(n-1)- \Ein([g])}< \min\bigl\{n-p-1,\frac{(n-2)^2}{n}\bigr\}.$$
To prove the last inequality, let $k=\min\bigl\{n-p-1,\frac{(n-2)^2}{n}\bigr\}$. The inequality $(n-2)\frac{n-1- \Ein([g])}{2(n-1)- \Ein([g])}<k$ is equivalent to $\Ein[g]>(n-1)\frac{n-2k-2}{n-k-2}$. Letting $k=n-p-1$, the inequality reads $\Ein[g]>\frac{(n-1)(2p-n)}{p-1}$. For $k=\frac{(n-2)^2}{n}$, the inequality reads $\Ein[g]>-\frac{(n-1)(n-4)}{2}$.\\
The corollary follows then directly from our Theorem B and Theorem 4.6(i) in  \cite{Schoen-Yau}. In fact, the later asserts that 
\[d<\min\bigl\{n-p-1,\frac{(n-2)^2}{n}\bigr\}\implies \pi_2(M)=...=\pi_p(M)=0.\]
This completes the proof of the corollary.
\end{proof}
Analogous results for the small $\ein$ invariant holds. We prove for instance the following one
\begin{proposition}
Let $(M,g)$ be a compact connected locally conformally flat Riemannian manifold of positive scalar curvature of dimension $n\geq 4$. Then 
$$ \ein([g])\leq \frac{2d-n+2}{d}.$$
\end{proposition}
\begin{proof}
We proceed as in the proof of Theorem B using Nayatani metric $\bar{g}$.
Recall that 
\begin{equation*}
\Eink(\bar{g})=\left((n-1)(n-2-2\delta)-k(n-2-\delta)\right)(\tr_{\bar{g}}{\mathcal A}) \bar{g}+k(n - 2)(\delta + 1){\mathcal A}.
\end{equation*}
The tensor ${\mathcal A}$ is symmetric and non negative, denote by $ {\mathcal A}_{\rm {max}}$ its maximum eigenvalue. This is a positive function as the trace of ${\mathcal A}$ is positive. Then it is easy to see that for $k<0$, the tensor $\Eink(\bar{g})$ is positive if 
$$\left((n-1)(n-2-2\delta)-k(n-2-\delta)\right)  {\mathcal A}_{\rm {max}}  +k(n - 2)(\delta + 1){\mathcal A}_{\rm {max}}>0.$$
Consequently, $\Eink(\bar{g})>$ for all $k<0$ satisfying $k>\frac{2\delta-n+2}{\delta}.$ This completes the proof of the proposition.
\end{proof}
We remark that the inequality in Proposition B is not optimal as one can check it for the case of the standard product $S^{n-d-1}\times H^{d+1}$. The author believes that one can improve this inequality to optimality using a deeper look at  the tensor ${\mathcal A}$.
\subsubsection{Proof of Corollary ${\rm B}^{\prime}$}
\begin{proof}
First we remark that 
$$\Ein([g])>(n-1)\frac{n-4}{n-3}\iff  (n-2)\frac{n-1-\Ein([g])}{2(n-1)-\Ein([g])}<1.$$
Consequently, under our assumption, Theorem B implies that $d(M,[g])<1$. The corollary follows then directly from Theorem 6.1 of Izeki \cite{Izeki}.
\end{proof}

\subsection{Connected sums of conformally flat manifolds of positive $\Eink$ tensor: Proof of Theorem C and Corollary C}

\subsubsection{Proof of Theorem C}

\begin{proof} This theorem follows from a general theorem due to Hoelzel, see Theorem 6.1 in \cite{Hoelzel}. We shall use the same notations as in \cite{Hoelzel}. Let $C_B({\Bbb{R}}^n)$ denote the  vector space of algebraic curvature operators $\Lambda^2 {\Bbb{R}}^n\rightarrow \Lambda^2 {\Bbb{R}}^n$ satisfying the first Bianchi identity and endowed with the canonical inner product.
Let 
$$C_{\Eink >0}:=\{R\in C_B({\Bbb{R}}^n): \Eink(R)>0\},$$
where for a unit vector $u$, $\Eink(R)(u)=\Scal(R) u-k\Ric(u)$. Here $\Ric$ and  $\Scal(R)$ denote respectively  as usual the first Ricci contracttion  and the full contraction of $R$. The subset $C_{\Eink >0}$ is clearly open, convex, star shaped with respect to the origin and it is an ${\rm{O}}(n)$-invariant cone. Furthermore, it is easy to check that $\Eink(S^{n-1}\times \Bbb{R})>0$ for $n-1>k$. The theorem follows then from Theorem 6.1 in \cite{Hoelzel}.

\end{proof}
\subsubsection{Proof of Corollary C}
\begin{proof}
Recall that the standard metric  on $S^1\times S^{n-1}$ is locally conformally flat and has its $\Ein$ and $\ein$  respectively equal to $n-1$ and $-\infty$. The corollary follows then from Theorem C.
\end{proof}

\section{$\Ein([g])$ and $\ein([g])$ invariants for general conformal classes in four dimensions: Proof of Theorems D and ${\rm D}^{\prime}$}
\subsection{Proof of Theorem D}

\begin{proof}
The first part is a direct consequence of Corollary B to Theorem A in \cite{CGY}. In fact, the condition $\frac{1}{2}\Scal. g-\Ric>0$ is equivalent to $\Ein_2(g)>0$ and therefore implies $\Ein(g)>2$. The remaining two parts are a consequence of the main Theorem in \cite{Gursky-Viacklovsky} which asserts that under the assumption
$$4\int_M\sigma_2(A)\mu_g+\frac{\alpha(\alpha+1)}{6}\left(Y[g]\right)^2>0,$$
for an arbitrary positive constant $\alpha$, the existence of a conformal metric $\bar{g}\in [g]$ with
\begin{equation}\label{GV}
\Ein_{\frac{-2}{\alpha}}(\bar{g})>0\quad {\rm and}\quad \Ein_{\frac{2}{\alpha+1}}(\bar{g})>0.
\end{equation}
So if  $\int_M\sigma_2(A)=0$ and since $\alpha>0$ is arbitrary we immediately get the conclusions $\Ein([g])\geq 2$ and $\ein([g])=-\infty$.
Finally, if $\int_M\sigma_2(A)<0$ denote by $c$ the following positive constant
$$c=\frac{-24\int_m\sigma_2(A)\mu_g}{\left(Y[g]\right)^2}.$$
The aforementioned theorem guarantees under the condition $\alpha>\frac{1}{2}\left( \sqrt{4c+1}-1\right)$ the existence of a conformal metric $\bar{g}\in [g]$ that has the property \ref{GV}, that is 
$\ein([g])<\frac{-2}{\alpha}$ and $ \Ein([g]) >\frac{2}{\alpha+1}.$
In other words, 
$$\alpha>\frac{1}{2}\left( \sqrt{4c+1}-1\right)\implies \alpha >\frac{-2}{\ein([g])}$$
and
$$\alpha>\frac{1}{2}\left( \sqrt{4c+1}-1\right)\implies \alpha > \frac{2}{\Ein([g])}-1.$$
From which we conclude that 
$$\ein([g])\leq \frac{-4}{\sqrt{4c+1}-1}\quad {\rm and}\quad \Ein([g])\geq \frac{4}{\sqrt{4c+1}+1}.$$
This completes the proof.
\end{proof}

\subsection{The positivity of Paneitz operator and the condition $\Ein([g])>1$}
\subsubsection{Proof of Theorem ${\rm D}^{\prime}$}
\begin{proof}
Note that the condition $\Ein([g])>1$ implies by definition the existence of a metric $g_1$ in the conformal class $[g]$ such that $\Eink(g_1)>0$ for some $k>1$ and the theorem follows from Theorem 5.5 in \cite{ES}. Also, since $\Eink(g_1)>0\implies {\rm Ein}_1(g_1)>0$ the theorem follows from  Proposition 6.1 in \cite{Gursky-Viacklovsky}. For the seek of completeness we provide a  proof of this theorem\\
Let $\lambda$ be an eigenvalue of $P_g$ and $u$ a corresponding eigenfunction. Using the fact that the Laplacian operator is self adjoint we get 
\begin{equation}
\begin{split}
\int <\lambda u,u>\mu_g=&\int<Pu,u>\mu_g=\lambda \int u^2\mu_g\\
=& \int\Big\{<\Delta^2u,u>+\bigl(\frac{2}{3}{\rm Scal}\, g-2 \Ric \bigr )(\nabla u,\nabla u)\Big\}\mu_g\\
=& \int\Big\{(\Delta u)^2+\bigl(\frac{2}{3}{\rm Scal}\, g-2 \Ric \bigr)(\nabla u,\nabla u)\Big\}\mu_g\\
=& \int\Big\{\frac{-1}{3}(\Delta u)^2+\frac{4}{3}(\Delta u)^2+\bigl(\frac{2}{3}{\rm Scal}\, g-2 \Ric \bigr)(\nabla u,\nabla u)\Big\}\mu_g.
\end{split}
\end{equation}
The Bochner formula shows that

\begin{equation}
\frac{4}{3}\int(\Delta u)^2\mu_g=\frac{4}{3}\int \big \{ |{\rm Hess}\, u|^2+\Ric (\nabla u,\nabla u)\big\}\mu_g.
\end{equation}
Hence, we deduce that

\begin{equation}
\begin{split}
\lambda \int u^2\mu_g=& \int\Big\{\frac{-1}{3}(\Delta u)^2+ \frac{4}{3}|{\rm Hess}\, u|^2+\bigl(\frac{2}{3}{\rm Scal}\, g-\frac{2}{3} \Ric \bigr)(\nabla u,\nabla u)\Big\}\mu_g\\
\geq & \frac{2}{3}\int \Big\{ \bigl({\rm Ein}_1 \bigr)(\nabla u,\nabla u)\Big\}\mu_g.
\end{split}
\end{equation}
Where in the last step we used Newton-Maclaurin's identity as follows
$$\left(\sigma_1({\rm Hess}\,u)\right)^2\geq 8/3\sigma_2({\rm Hess}\,u)= 4/3\left( \left(\sigma_1({\rm Hess}\,u)\right)^2-|{\rm Hess}\, u|^2\right).$$
 The first part of the result follows directly. If $\lambda =0$ then clearly $\nabla u=0$ and then the function $u$ is constant.
\end{proof}

\subsubsection{Proof of Corollary ${\rm D}^{\prime}$}
\begin{proof}
First note that the previous proposition guarantees that the kernel of $P_g$ consists of constant functions.\\
Next, since the condition  $\Ein([g])>1$ implies the positivity of the Yamabe constant of the conformal class $[g]$,   Theorem B of Gursky \cite{Gursky} shows that 
$\int Q_g\mu_g \leq 8\pi^2$ with equality if and only if the manifold is conformally equivalent to the sphere. In the case $\int Q_g\mu_g < 8\pi^2$,  
a theorem of Djadli-Malchiodi \cite{DjMa} guarantees the existence of a metric with constant $Q$-curvature.
\end{proof}
\begin{remarks}
\begin{enumerate}
\item M. Lai proved in \cite{Lai} an equivalent version the above theorem and corollary under the condition of $3$-positive Ricci curvature. It is easy to show that  in four dimensions, a metric $g$ has $3$-positive Ricci curvature if and only if ${\rm Ein}_1(g)>0$.
\item For higher dimensions $n\geq 4$, let $k=\frac{2n(n-1)}{3n-4}$. It is proved in \cite{modified-Eink} that $\Eink>0\implies \Gamma_2(A)>0$, that is positive scalar curvature and positive $\sigma_2$ curvature, in particular it implies the positivity of the integral of the $Q$-curvature.
\end{enumerate}
\end{remarks}


\begin{thebibliography}{9}
\bibitem{Besse} Besse A. L., Einstein Manifolds, Springer, Berlin-New York (1987).
%\bibitem{Bishop} Richard L. Bishop, Richard J. Crittenden, Geometry of Manifolds, Academic Press, 1964.
%\bibitem{Bolotov} Bolotov, D.,   Gromov’s macroscopic dimension conjecture, Algebraic and Geometric Topology 6 (2006) 1669–1676.
%\bibitem{BoLa} Botvinick B. and Labbi M. L., Compact manifolds with positive $\Gamma_2$ curvature, Journal of Differential Geometry and its Applications, ...
%\bibitem{Bourguignon} J.-P. Bourguignon, Ricci curvature and Einstein metrics, Global differential geometry and global analysis (Berlin, 1979), Lecture Notes in Math., vol. 838, Springer,Berlin, 1981, pp. 42–63
%\bibitem{DjMa} Djadli Z.  and  Malchiodi A., Existence of conformal metrics with constant Q-curvature,
%Ann. of Math. (2) 168 (2008), no. 3, 813–858.
%\bibitem{ES} Eastwood M. G. and Singer M., The Frohlicher spectral sequence on a twistor space,  Journal of Differential Geometry, 38 (1993) 653 669.

\bibitem{Bourguignon} Bourguignon J. P., Ricci curvature and Einstein metrics, Global differential geometry and global analysis (Berlin, 1979), Lecture Notes in Math., vol. 838, Springer, Berlin,  42-63, (1981).


\bibitem{CGY} S.-Y. A. Chang, M. J. Gursky and  P. C. Yang, An equation of Monge-Ampere type in
conformal geometry, and four-manifolds of positive Ricci curvature, Ann. of Math.
(2) 155(3) (2002) 709-787.
%\bibitem{Chang-Hang-Yang} S.-Y. A. Chang, F. Hang, and P. C. Yang., On a class of locally conformally flat
%manifolds,  Int. Math. Res. Not., (4) 185-209, (2004).

%\bibitem{CTZ} Chen B-L., Tang S-H. and Zhu X-P.,  Complete classification of compact four manifolds with positive isotropic curvature, J. %Differential Geometry, 91 (2012), 41-80.
\bibitem{DjMa} Djadli Z.  and  Malchiodi A., Existence of conformal metrics with constant Q-curvature,
Ann. of Math. (2) 168 (2008), no. 3, 813-858.
\bibitem{ES} Eastwood M. G. and Singer M., The Frohlicher spectral sequence on a twistor space,  Journal of Differential Geometry, 38 (1993) 
653-669.
%\bibitem{Gonzalez} Gonzalez M. d. M., Singular sets of a class of locally conformally
%flat manifolds, Duke Math. J. 129(2005), no. 3, 551–572.
%\bibitem{GroLa} M. Gromov and B. Lawson, classification of simply connected manifolds of positive scalar curvature
%\bibitem{Gromov} M. Gromov, A Dozen Problems, Questions andConjectures about Positive Scalar Curvature, preprint.
\bibitem{Guan-Wang} Guan P.,  Wang G., Conformal Deformations of the Smallest Eigenvalue of the Ricci Tensor,  American Journal of Mathematics,
Vol. 129, No. 2 (2007), 499-526.
\bibitem{Gursky} Gursky M., The Principal Eigenvalue of a Conformally Invariant Differential Operator,
with an Application to Semilinear Elliptic PDE, Commun. Math. Phys. 207, 131-143 (1999).
\bibitem{Gursky-Viacklovsky} Gursky M.,  Viacklovsky J. A., A fully nonlinear equation on four-manifolds with positive scalar curvature, J. differential geometry 63 (2003) 131-154
\bibitem{Hoelzel} Hoelzel S., Surgery stable curvature conditions, Math. Ann. 365, 13-47 (2016). 
\bibitem{Izeki} Izeki H., Limit sets of Kleinian groups and conformally flat Riemannian manifolds, Invent. math., 122, 603-625 (1995).
\bibitem{Izeki2} Izeki H., Convex-cocompactness of Klenian groups and conformally flat manifolds with positive scalar curvature, Proc. Am. Math. Soc. 130, no. 12, 3731-3740 (2002).

\bibitem{modified-Eink} Labbi M., On modified Einstein tensors and two smooth invariants of compact manifolds, Transactions of the American Mathematical Society (accepted for publication).
%\bibitem{agag} Labbi M.L., Stability of the p-curvature positivity under surgeries and manifolds with positive Einstein tensor, Annals of Global analysis and geometry, 15: 299-312, 1997.
\bibitem{Labbi-Betti}  Labbi M., Sur les nombres de Betti des vari\'et\'es conform\'ement plates,
Comptes rendus de l'Acad\'emie des sciences. S\'erie 1, Math\'ematique 319 (1), 77-80 (1994).

%\bibitem{Labbi-einstein} M.-L Labbi, Compact manifolds with positive
       % Einstein curvature. Geom. Dedicata 108 (2004), 205-217.

\bibitem{Labbi-nachr} Labbi M., On Weitzenbock curvature operators, Mathematische Nachrichten 288 (4), 402-411 (2015).

\bibitem{Lafontaine} Lafontaine J., Remarques sur les vari\'et\'es conform\'ement plates, Math. Ann. 259 (1982), 313–319.

\bibitem{Lai}  Lai M., A remark on the nonnegativity of the Paneitz operator,  Proc. Amer. Math. Soc. 143 (2015), 4893-4900.

\bibitem{Nayatani} Nayatani S., Patterson-Sullivan measure
and conformally flat metrics, Math. Z. 225, 115-131 (1997).

\bibitem{Noronha} Noronha M. H., Some compact conformally flat manifolds with nonnegative scalar curvature, Geometria Dedicata, 47, 255-268 (1993).

\bibitem{Rastall} Rastall P., Generalization of the Einstein theory, 
Phys. Rev. D, 6, 3357-3359, (1972).

\bibitem{Schoen-Yau}  Schoen R. and Yau S. T., Conformally flat manifolds, Kleinian groups
and scalar curvature, Invent. math. 92, 47-71 (1988).
\bibitem{Schoen-Yau2} Schoen R. and Yau S. T., Lectures on differential geometry.  Conference Proceedings
and Lecture Notes in Geometry and Topology, I. International Press, 1994.



\end{thebibliography}
\end{document}